\documentclass[amscd,syntonly,amssymb]{amsart}

\usepackage{epsfig}
\theoremstyle{plain}
\newtheorem{Thm}{Theorem}
\newtheorem{Prop}[Thm]{Proposition}
\newtheorem{Cor}[Thm]{Corollary}
\newtheorem{Lem}[Thm]{Lemma}

 \theoremstyle{definition}
\theoremstyle{remark}

\numberwithin{equation}{section}

\begin{document}
 \title{A Characteristic Number of Hamiltonian Bundles over $S^2$}

 \author{ ANDR\'{E}S   VI\~{N}A}
\address{Departamento de F\'{i}sica. Universidad de Oviedo.   Avda Calvo
 Sotelo.     33007 Oviedo. Spain. }
 \email{vina@uniovi.es}
\thanks{This work has been partially supported by Ministerio de Ciencia y
Tecnolog\'{\i}a, grant MAT2003-09243-C02-00}

 \maketitle
\begin{abstract}
Each loop $\psi$ in the group $\text{Ham}(M)$ of Hamiltonian
diffeomorphisms  of a symplectic manifold $M$ determines a
fibration $E$ on $S^2$, whose coupling class \cite{G-L-S} is
denoted by $c$. If $VTE$ is the vertical tangent bundle of $E$, we
relate the characteristic number $\int_E c_1(VTE)c^n$ with the
Maslov index of the linearized flow $\psi_{t*}$ and the Chern
class $c_1(TM)$. We give the value of this characteristic number
for loops of Hamiltonian symplectomorphisms of Hirzebruch
surfaces.

\end{abstract}
\section {Introduction} \label{S:intro}

A  loop $\psi:S^1\to\text{Ham}(M,\omega)$ in the group of
Hamiltonian diffeomorphisms \cite{Mc-S} of a symplectic manifold
$(M^{2n},\omega)$ can be considered as a clutching function of  a
Hamiltonian fibration $E\stackrel{\pi}{\to}S^2$ with fibre $M$.
The total space $E$ supports the coupling class $c\in H^2(E,{\Bbb
R})$; this is the unique class such that $c^{n+1}=0$, and
$i_p^*(c)$ is the cohomology class of the symplectic structure on
the fibre $\pi^{-1}(p)$, where
 $i_p$ is the inclusion of  $\pi^{-1}(p)$ in $E$ \cite{G-L-S}.
Furthermore  one can consider on $E$ the first Chern class
$c_1(VTE)$ of the vertical tangent bundle of $E$. These canonical
cohomology classes on $E$ determine the characteristic number (see
\cite{L-M-P})
\begin{equation}\label{defIpsi}
I_{\psi}=\int_E c_1(VTE)\,c^n,
\end{equation}
 which depends only on the homotopy class of $\psi$.
   Since $I$ is an ${\Bbb R}$-valued group homomorphism on
$\pi_1(\text{Ham}(M,\omega))$, the non vanishing of $I$ implies
that the   group $\pi_1(\text{Ham}(M,\omega))$ is infinite. That
is, $I$ can be used to detect the infinitude  of the corresponding
homotopy group. Furthermore $I$ calibrates the Hofer's norm $\nu$
on $\pi_1(\text{Ham}(M,\omega))$ in the sense that $\nu(\psi)\geq
C|I_{\psi}|$, for all $\psi$, where $C$ is a positive constant
\cite{lP01}.

$I$ is a generalization of the mixed action-Maslov homomorphism
introduced by Polterovich \cite{lP97} for {\it monotone}
manifolds, that is, when $[\omega]=a c_1(TM)$ and $a
>0$. The value of this mixed action-Maslov homomorphism on a loop $\psi$ is, in
many cases, easy to calculate, since it is a linear combination of
the symplectic action
 around {\it any} orbit
$\{\psi_t(x_0)\}_t$ and the Maslov index of the linearized flow
$(\psi_t)_*$ along this orbit.
 By contrast, $I$ is defined for
Hamiltonian loops in general manifolds (non necessarily monotone),
and  its value is mostly not so easy to determine from the
definition.

Our purpose in this note is to obtain an explicit expression for
$I_{\psi}$, which can be used to calculate its value.
 More precisely,
 when the bundle $TM$ admits local symplectic
trivializations whose  domains are fixed by the diffeomorphisms
$\psi_t$, we deduce a formula for $I_{\psi}$ in which appear a
contribution related to the Maslov indices  of the linearized flow
$\psi_{t*}$ in the trivializations, and a
 second one in which are involved  transition functions of the
 bundle $\text{det}(TM)$.
 The second contribution is related with
the Chern class $c_1(TM)$ in the following sense. Using the
expression of $c_1(M)$ in terms of the transition functions of
$TM$ determined by the trivializations, $\langle
c_1(M)[\omega]^{n-1},\,M\rangle$ can be written as a sum
$\sum_j\int_{R_j}\sigma_j$, where $\sigma_j$ is a $2n-1$ form (see
(\ref{lasR})).  It turns out that the second contribution is equal
to this sum ``weighted" by a multiple of the Hamiltonian $f_t$
which generates $\psi$; more concretely, that contribution is
  $-n\sum_j \int dt\int_{R_j}(f_t\circ\psi_t)\sigma_j$.

Let $(M,\omega,f)$ be an integrable system such that the points
where the integrals of motion are dependent form a set $P$ which
is union of codimension $2$ submanifolds of $M$, and such that
$M\setminus P$ is invariant under $\psi_t$ and on it there exist
action-angle coordinates. Furthermore we assume that there are
$\psi_t$ invariant Darboux charts which cover $P$. Then
 the expression of
$I_{\psi}$ in this atlas reduces to the aforesaid second
contribution; that is, $I_{\psi}=-n\sum_j \int_{R_j}f\sigma_j$.

 The paper is organized as follows. In Section 2 we recall the
construction of the coupling class $c$ following \cite{lP01}.
Section 3 is concerned with the proof of the mentioned expression
for $I_{\psi}$. First we express $\langle
c_1(M)[\omega]^{n-1},\,M\rangle$ as the sum
$\sum_j\int_{R_j}\sigma_j$ of integrals of $2n-1$ forms, and next
we use this result to prove  the  formula for the invariant
$I_{\psi}$.
 In Section 4 we check and apply the formulae obtained in
Section 3. Using these formulae, we calculate $I_{\psi}$, when
$\psi$ is
 the loop in $\text{Ham}(S^2)$  generated by the $1$-turn rotation of $S^2$
around the $z$-axis. The result $I_{\psi}=0$ agrees with the fact
that $\pi_1(\text{Ham}(S^2))={\Bbb Z}_2$ and $I$ is a group
homomorphism on $\text{Ham}(M)$. We also prove that $I$ on
$\pi_1(\text{Ham}({\Bbb T}^{2n}))$ vanishes identically.
 When  $n=1$ this result is consistent with the
fact that $\pi_1(\text{Ham}({\Bbb T}^{2}))=0$. Finally we
determine the value of $I$ on the loops generated by action of
${\Bbb T}^2$ on  a general symplectic Hirzebruch surface (see
Theorem \ref{ThemHiz}).

 I thank Dusa McDuff for explaining me properties of the Maslov index of the linearized flow, and
 Eva Miranda  for clarifying me some points relative to  action-angle variables.


  \smallskip

\section{The coupling class}\label{S:Repre}

Let $(M,\omega)$ be a compact connected symplectic $2n$-manifold.
Let $\psi:S^1={\Bbb R}/{\Bbb Z}\to \text{Ham}(M,\omega)$ be a loop
in the group $\text{Ham}(M,\omega)$ at $\text{id}$. By $X_t$ is
denoted the  time-dependent vector field generated by $\psi_t$ and
$f_t$ is the normalized
 time-dependent Hamiltonian; that is,
$$\frac{d\psi_t}{dt}=X_t\circ\psi_t,\;\;\;\iota_{X_t}\omega=-df_t,\;\;\;\int_M f_t\omega^n=0.$$

Given $\epsilon$, with $0<\epsilon<\pi/2$, we set
$$D^2_+:=\{p\in S^2\,|\,0\leq\theta(p)<\pi/2+\epsilon \}$$
$$D^2_-:=\{p\in S^2\,|\,\pi/2-\epsilon<\theta(p)\leq\pi \},$$
 where $\theta\in[0,\,\pi]$ is the
  polar angle from the z-axis.

  Next we construct the Hamiltonian bundle $E$ over $S^2$
determined by $\psi$. First of all we extend $\psi$ to a map
defined on $F:=D^2_+\cap D^2_-$ by putting
$\psi(\theta,\phi)=\psi_t$, with $t=\phi/2\pi$, with $\phi$ the
spherical azimuth angle. We set
$$E=[(D^2_+\times M)\cup(D^2_-\times M)]/\simeq,\;\;\
\text{where}$$
$$(+,p\,,x)\simeq(-,p',y)\;\;\text{iff}\;\;
\begin{cases} p=p'\in F,\\
 y=\psi_t^{-1}(x),\;t=\phi(p)/2\pi.
 \end{cases}$$
In this way $M\hookrightarrow E\stackrel{\pi}{\rightarrow}S^2$ is
 a Hamiltonian bundle over $S^2$.

  We assume that $D^2_{\pm}$ are endowed
with the orientations induced by the usual one of $S^2$ (that is,
the orientation of $S^2$ as border of the unit ball). We suppose
that $S^1$ is oriented by $dt=d\phi/2\pi$, that is, $S^1$ is
oriented as $\partial D_{+}$ . In $E$ one considers the
orientation induced by the one defined on $D^2_{+}\times M$ by
$d\theta\wedge d\phi\wedge \omega^n$.

Let $\alpha$ be a monotone smooth map
$\alpha:[\pi/2-\epsilon,\,\pi]\to[0,\,1]$, with $\alpha(\theta)=1$
for $\theta\in[\pi/2-\epsilon,\,\pi/2+\epsilon]$ and
$\alpha(\theta)=0$ for $\theta$ near $\pi$. Now we consider the
$2$-form (see \cite{lP01})
 \begin{equation}\label{taudf}
 \tau=\begin{cases}
\omega,\;\;\text{on}\;\; D^2_+\times M \\
\omega+d(\alpha(f_t\circ\psi_t))\wedge dt,\;\; \text{on}\;\;
D^2_-\times M.
\end{cases}
\end{equation}
As $\alpha$ vanishes near $\pi$, $\tau$ is well-defined on
$D_{-}^2\times M$; moreover  on $F\times M\subset D^2_-\times M$,
$\tau$ reduces to $\omega+d(f_t\circ\psi_t)\wedge dt$. If we
denote by $h$ the
 map
 $$h:F\times M\subset D^2_-\times M\rightarrow F\times M\subset D^2_+\times
 M$$
given by $h(p,x)=(p,\psi_t(x))$, with $t=\phi(p)/2\pi$, then
taking into account that $h_*(\frac{\partial}{\partial
t})=\frac{\partial}{\partial t}+X_t\circ\psi_t$, it follows from
$\iota_{X_t}\omega=-df_t$ that
 $h^*\omega= \omega+d(f_t\circ\psi_t)\wedge dt$. So one has the following Proposition
\begin{Prop}$\tau$ defines a closed $2$-form on $E$.
\end{Prop}

Moreover the cohomology class $[\tau]\in H^2(E,\,{\Bbb R})$ restricted to each fibre coincides with $[\omega]$.
On the other hand
$$\int_E\tau^{n+1}=
(n+1)\int_{D^2_{-}\times
M}(f_t\circ\psi_t)\alpha'(\theta)d\theta\wedge
dt\wedge\omega^{n}.$$ From the normalization condition for $f_t$
it follows that $\int_E\tau^{n+1}=0$. Hence $[\tau]$ is the
coupling class $c$ of the fibration $E$ \cite{G-L-S} \cite{Mc-S}.


\section{The characteristic number $I_{\psi}$.}
  Denoting $TM=\{v_x\in T_x M\,|\,x\in
M\}$, we put
$$VTE=[(D^2_+\times
TM)\cup(D^2_-\times TM  )]/\simeq,$$
 with
$$(+,p\,,v_x)\simeq(-,p',v'_{x'})\;\text{iff}\;
p=p',\;x'=\psi_t^{-1}(x),\; v'_{x'}=(\psi_t^{-1})_*(v_x)$$ where
$t=\phi(p)/2\pi$. So $VTE$ is a vector bundle over $E$; by
construction it is the vertical tangent bundle of $E$.

Let $(U; X_1,\dots,X_{2n})$ be a symplectic trivialization of $TM$
on $U\subset M$, and $(V; Y_1,\dots,Y_{2n})$ be a symplectic
trivialization on $V\subset M$. We put
 \begin{equation}\label{Upm}
 U_{\pm}:=\{
[\pm,p,x]\,|\,p\in D^2_{\pm}, x\in U \}
 \end{equation}
 and similarly for $V_{\pm}$.Denoting $x_t:=\psi_t^{-1}(x)$ one
has
$$U_+\cap U_-=\{[+,p,x]\,|\,p\in F,\,x\in U,\,x_t\in U  \}$$
$$V_+\cap V_-=\{[+,p,x]\,|\,p\in F,\,x\in V,\,x_t\in V  \}$$
$$V_-\cap U_-=\{[-,p,x]\,|\,p\in D^2_-,\,x\in V\cap U \}$$
$$U_+\cap V_+=\{[+,p,x]\,|\,p\in D^2_+,\,x\in V\cap U \}$$
The corresponding transition functions of $VTE$ are
$$g_{U_-U_+}([+,p,x])=A(t,x)\in Sp(2n,{\Bbb R}),\; \text{with}\;
 \psi_{t*}^{-1}\big(X_i(x)\big)=
 \sum_k A^k\,_i(t,x)\,X_k(x_t)$$
$$g_{V_-V_+}([+,p,x])=B(t,x)\in Sp(2n,{\Bbb R}),\; \text{with}\;
 \psi_{t*}^{-1}\big(Y_i(x)  \big)=
 \sum_k B^k\,_i(t,x)\, Y_k(x_t)$$
 $$g_{U_-V_-}([-,p,x])=R(x)=g_{U_+V_+}([+,p,x]),\; \text{with}\;
 Y_i(x)= \sum_k R^k\,_i(x)\,X_k(x).$$

 We denote by $\rho$ the usual map $\rho:Sp(2n,\,{\Bbb R})\to U(1)$
 which restricts to the determinant map on $U(n)$
 \cite{S-Z}, then $l_{ab}:=\rho\circ g_{ab}$ is a transition function
 for $\text{det}(VTE)$. We also use the following notation, the
 matrices in $Sp(2n,{\Bbb R})$ are denoted with capital letters and its
 images by $\rho$ will be denoted by the corresponding small
 letters; that is,
 \begin{equation}\label{a-b-r}
 a(t,x):=\rho(A(t,x)),\;\;\; b(t,x):=\rho(B(t,x)) ,\;\;\;
 r_{UV}(x):=\rho(R(x)).
 \end{equation}

If $\psi_t(U)\subset U$ for all $t$, given $x\in U$, the winding
number of the map $t\in S^1\mapsto a^{-1}(t,x)\in U(1)$ is the
integer
 \begin{equation} \label{Judef}\frac{i}{2\pi }\int_{0}^1 a^{-1}(t,x)\frac{\partial
a}{\partial t}(t,x) dt.
 \end{equation}
 This integer is independent of the point
$x\in U$, it will be denoted $J_U$. The number $J_U$ is the Maslov
index in $U$ of the linearized flow $\psi_{t*}$.
 Analogously, if $\psi_t(V)\subset V$ for all $t$   we have the integer
\begin{equation} \label{Jvdef}
J_V=\frac{i}{2\pi }\int_{0}^1b^{-1}(t,x)\frac{\partial b}{\partial
t}(t,x) dt,
 \end{equation}
 $x$ being any point of $V$; this is the
Maslov index in $V$ of  $\psi_{t*}$.

As a previous step to compute $I_{\psi}$ we shall prove the
following Lemma, in which the value $\langle
c_1(M)[\omega]^{n-1},\, [M]\rangle$ is expressed in terms of
transition functions of $\text{det}(TM)$.

 \begin{Lem}\label{Lemaux}
 Let $\{B_1,\dots, B_m  \}$ be a set of trivializations of $TM$,
 such that its domains cover $M$. Then
 \begin{equation}\label{Lem0}\langle
 c_1(TM)[\omega]^{n-1},\,[M]\rangle=\frac{-i}{2\pi}\sum_{i<k}\int_{A_{ik}}d (\log
 s_{ik})\wedge \omega^{n-1},
 \end{equation}
 $s_{ik}$ being the corresponding transition function of
 $\text{det}\,(TM)$ and
 \begin{equation}\label{Aik}
 A_{ik}=\Big(\partial B_i\setminus\cup_{r<k}B_r  \Big)\cap B_k.
  \end{equation}
\end{Lem}
\begin{proof}
$c_1(M)$ is represented on $B_a$ by the $2$-form
$$\frac{-i}{2\pi}\sum_c d\big(\varphi_cd\,\log s_{ac}\big),$$
where $\{ \varphi_c\}$ is a partition of unity subordinate to the
covering $\{B_1,\dots,B_m\}$.

If $m=2$
$$\langle
 c_1(M)[\omega]^{n-1},\,[M]\rangle=
 \frac{-i}{2\pi}\int_{B_1}d\big(\varphi_2 d\,\log
 s_{12}\big)\wedge \omega^{n-1}+
\frac{-i}{2\pi}\int_{B_2\setminus B_1}d\big(\varphi_1 d\,\log
 s_{21}\big)\wedge \omega^{n-1}.$$
By Stokes' theorem
 \begin{equation}\label{Lm1} \langle c_1(M)[\omega]^{n-1},\,[M]\rangle=\int_{\partial
 B_1}\varphi_2L_{12}+ \int_{\partial
 ({B_2\setminus B_1})}\varphi_1L_{21},
 \end{equation}
 where
 $$L_{jk}:=(-i/2\pi)d\,\log s_{jk}\wedge \omega^{n-1}.$$
 Since $\partial
 ({B_2\setminus B_1})\cap B_1=\emptyset$, $\varphi_1$ vanishes on
$\partial
 ({B_2\setminus B_1})$ and the last integral in (\ref{Lm1}) is
 zero.

As $\varphi_2$ is $1$ on $\partial B_1$, we have
$$\langle c_1(M)[\omega]^{n-1},\,[M]\rangle=\int_{\partial
 B_1}L_{12}.$$
 In this case $\partial B_1\subset B_2$, so $\partial B_1=A_{12}$,
 and the the Lemma is proved when $m=2$.

 If $m=3$
\begin{align}\label{star}
\langle c_1(M)[\omega]^{n-1},\,[M]\rangle &=  \int_{\partial
 B_1}\big(\varphi_2L_{12}+\varphi_3L_{13}  \big)+ \int_{\partial
 ({B_2\setminus B_1})}\big( \varphi_1L_{21}   +\varphi_3L_{23}\big) \\ \label{starsig} &+
 \int_{\partial
 ({B_3\setminus (B_1\cup B_2)})}\big(\varphi_1L_{31}+\varphi_2L_{32}\big).
  \end{align}
  As $ \partial
 ({B_3\setminus (B_1\cup B_2)})$ and the interior of $B_1\cup B_2$ are disjoint
 sets, $\varphi_1$ and $\varphi_2$ vanish on
 $ \partial({B_3\setminus (B_1\cup B_2)})$, and the  integral
 in (\ref{starsig}) is zero. Analogously $\partial (B_2\setminus B_1)$
 and support of $\varphi_1$ are disjoint so
\begin{equation}\label{starbis} \langle
c_1(M)[\omega]^{n-1},\,[M]\rangle = \int_{\partial
 B_1}\big(\varphi_2L_{12}+\varphi_3L_{13}  \big)+ \int_{\partial
 ({B_2\setminus B_1})} \varphi_3L_{23}.
\end{equation}

On the other hand  $\partial B_1=A+D$, with $A:=\partial
B_1\setminus B_2$
 (oriented as $\partial B_1$) and $D:=(\partial B_1\setminus B_2)\cap
 B_3$ (see Figure 1).
 Moreover $\partial (B_2\setminus B_1)=-A+C$ with $C:=(\partial B_2\setminus B_1)\cap
 B_3$ (oriented as $\partial B_2$).

\begin{figure}[htbp]
\begin{center}
\epsfig{file=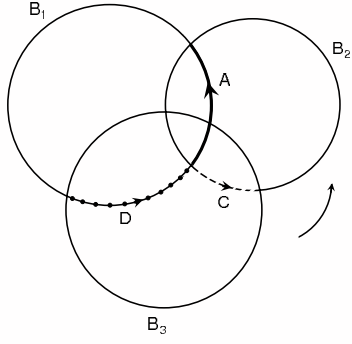,height=5cm}
\end{center}
\caption[Figure 1]{\small $A=\partial B_1\cap B_2$, $D=(\partial
B_1 \setminus B_2)\cap B_3$  and $C=(\partial B_2 \setminus
B_1)\cap B_3$.  }
\end{figure}

 Since $C\cap(B_1\cup B_2)=\emptyset$, then $\varphi_3|_C=1$; thus
\begin{equation}\label{aux1}
\langle c_1(M)[\omega]^{n-1},\,[M]\rangle =
\int_{A+D}\big(\varphi_2L_{12}+\varphi_3L_{13}  \big)+
\int_{-A}\varphi_3 L_{23} + \int_{A_{23}} L_{23}.
 \end{equation}
 The last integral in (\ref{aux1}) is just the term in (\ref{Lem0})
 with $i=2, k=3$.

Since $\varphi_j|_D=0$, for $j=1,2$, then $\varphi_3|_D=1$. As $A$
and support of $\varphi_1$ are disjoint sets, then
 $(\varphi_2+\varphi_3)|_A=1$. It follows
from these facts together with the cocycle condition
$L_{13}+L_{32}=L_{12}$ that
 \begin{equation}\label{Forl}\langle c_1(M)[\omega]^{n-1},\,[M]\rangle = \int_{A}L_{12} +
\int_{D}L_{13}+ \int_{A_{23}}L_{23}.
\end{equation}

On the other hand $A_{12}=(\partial B_1\setminus B_1)\cap B_2=A$.
Similarly $A_{13}=D$. Therefore (\ref{Forl}) is the formula given
in the statement of Lemma when $m=3$.

The preceding arguments can be generalized to any $m$
\begin{align}\label{Tocho}
\langle c_1(TM)[\omega]^{n-1},\,[M]\rangle & = \int_{\partial
B_1}\sum_{j\ne 1}\varphi_j L_{1j}+\dots+ \int_{\partial
(B_{m-1}\setminus \cup_{r<m-1}B_r)}\sum_{j\ne m-1}\varphi_j L_{m-1,j} \\
\label{caca} & +\int_{\partial (B_{m}\setminus
\cup_{r<m}B_r)}\sum_{j\ne m}\varphi_j L_{m-1,j}.
\end{align}
For any $j=1,\dots, m-1$ support of $\varphi_j$ and
$\partial\big(B_{m}\setminus \cup_{r<m}B_r\big)$ are disjoint
sets. Thus the integral (\ref{caca})   is zero (as in the cases
$m=2,3$). We decompose
$$\partial(B_{m-1}\setminus \cup_{r<m-1}B_r)=E+G,$$
with
$$E:= \Big(\partial B_{m-1}\setminus
\cup_{r<m-1}B_r\Big)\cap B_m.$$
    Then $\varphi_j|_E=0$ for all $j\ne m$ and $\varphi_m|_E=1$;
 thus
    \begin{equation}\label{tocho1}
    \int_{\partial
(B_{m-1}\setminus \cup_{r<m-1}B_r)}\sum_{j\ne m-1}\varphi_j
L_{m-1,j}=\int_G+\int_{A_{m-1,m}}L_{m-1,m}.
\end{equation}
    The last integral in (\ref{tocho1}) is the term in (\ref{Lem0})
    which corresponds to $i=m-1, k=m$. An analogous, but more tedious, calculation to
    the one for the case    $m=3$  allows to
 identify in (\ref{Tocho})
    the remainder
    terms of (\ref{Lem0}).

\end{proof}

Lemma \ref{Lemaux} gives a way for expressing $\langle
 c_1(TM)[\omega]^{n-1},\,[M]\rangle$ as a sum of integrals of   $2n-1$ differential forms on $2n-1$
 chains. The righthand side of (\ref{Lem0})  can be written
  schematically
\begin{equation}\label{lasR}
 \sum_j\int_{R_j}\sigma_j.
\end{equation}
In next Theorem we use this expression to give an explicit formula
for $I_{\psi}$ in terms of transition functions of
$\text{det}(TM)$ and Maslov indices of $\psi_{t*}$.

 \begin{Thm}\label{Thefl} If $\{B_1,\dots,B_m\}$ is a set of symplectic
trivializations for $TM$ which covers $M$, and such that
$\psi_t(B_j)=B_j$, for all $t$ and all $j$, then
\begin{equation}\label{subs}
I_{\psi}=\sum_{i=1}^m J_i\int_{B_i\setminus\cup_{j<i}B_j}
\omega^n+\sum_{i<k}N_{ik},
    \end{equation}
where
$$N_{ik}=n\frac{i}{2\pi}\int_{0}^{1}dt\int_{A_{ik}}(f_t\circ\psi_t)
(d\, \log r_{ik})\wedge \omega^{n-1},$$
 $A_{ik}=(\partial
B_i\setminus\cup_{r<k}B_r)\cap B_k$, $J_{i}$ is the Maslov index
of $(\psi_t)_*$ in the trivialization $B_i$ and $r_{ik}$ the
corresponding transition function of $\text{det}\,(TM)$.

\end{Thm}

\begin{proof}
 Using the notation (\ref{Upm}) we put
\begin{equation}\label{Opm}
O_{\bf{2a-1}}:=(B_{a})_{-},\;\;\,\; O_{\bf{2a}}:=(B_a)_{+}.
 \end{equation}
Then $\{O_{\bf{c}}\,|\, c=1,\dots,2m \}$ is a covering for $E$. We
shall denote by $l_{\bf{bc}}$ the respective transition functions
for $\text{det}\,(VTE)$. If we set $U:=B_1, \, V:=B_2$, one has by
(\ref{a-b-r})
$$l_{\bf{12}}=a(t,x),\;\; l_{\bf{13}}=r_{UV}(x)
,\;\;l_{\bf{34}}=b(t,x).$$
 We can determine $I_{\psi}=\langle c_1(VTE)\,c^n,\,[E]\rangle$
 applying the result given in Lemma \ref{Lemaux} to the set
 $\{O_{\bf{c}}\}$ of trivializations of $VTE$. That is,
 \begin{equation}\label{ITH}
 I_{\psi}=\sum_{\bf{a}<\bf{b}}\mathcal{
 T}_{\bf{ab}},\;\;\text{where}\;\;
\mathcal
 {T}_{\bf{ab}}=\frac{-i}{2\pi}\int_{A_{\bf{ab}}}d\,\log l_{\bf{ab}}\wedge
 \tau^n.
 \end{equation}

It follows from (\ref{Opm}) and (\ref{taudf}) that $\tau$ is equal
to $\omega$ on $A_{\bf{ab}}$ unless $\bf{a}$ and  $\bf{b}$ are
both odd; in this case
$\tau=\omega+d(\alpha(f_t\circ\psi_t))\wedge dt.$

We will calculate the summand  $\mathcal
 {T}_{\bf{12}}$ in (\ref{ITH}). The set $A_{\bf{12}}=\partial
O_{\bf{1}}\cap O_{\bf{2}}=\partial U_{-}\cap U_{+}$, and
$$\partial U_{-}=\{[+,p,x]\,|\,p\in\partial D^2_{-} ,\; x\in U \}\cup
\{[-,p,x]\,|\, p\in D^2_{-},\; x\in\partial U  \}.$$
 So
 $$A_{\bf{12}}=\{[+,p,x]\,|\,p\in\partial D^2_{-},\,x\in U  \}.$$
 Taking into account (\ref{Judef}), (\ref{a-b-r}) together with the fact that  orientations of $S^1$
 and $\partial D^2_{-}$ are opposite , we deduce
 $$ \mathcal{T}_{\bf{12}}=
 \frac{-i}{2\pi}\int_U\Big(\int_{-S^1}a^{-1}(t,x)\frac{\partial
a(t,x)}{\partial t} \,dt  \Big)\omega^n
  =J_U\int_U\omega^n.$$

 Next we consider the term  $\mathcal
 {T}_{\bf{34}}$. The integration domain is
 $$A_{\bf{34}}=\big(\partial V_{-}\setminus(U_{-}\cup U_{+}) \big)\cap
 V_{+}=\{[+,p,x]\,|\,p\in\partial D^2_{-},\; x\in V\setminus U  \}.$$
 Hence
$$\mathcal{T}_{\bf{34}} =J_V\int_{V\setminus U}\omega^n.$$

 In general,
 \begin{align} \notag A_{\bf{2j-1,2j}}&=\Big(\partial B_{j-}\setminus \cup_{r<j}\big(B_{r+}\cup B_{r-}\big)  \Big)\cap
 B_{j+}  \\ \notag &=\{[+,p,x]\,|\,p\in\partial D^2_-,\,x\in B_j\setminus\cup_{r<j}B_r  \}.
 \end{align}
  Hence the term in (\ref{ITH}) with $\bf{a}=2j-1$,
 $\bf{b}=2j$ gives a contribution to $I_{\psi}$ equal to
 \begin{equation}\label{2stars}
 J_{B_j}\int_{B_j\setminus U_{r<j}B_r}\omega^n
 \end{equation}

Now we analyze $\mathcal {T}_{\bf{13}}$.
 $$A_{\bf{13}}=\{[-,p,x]\,|\, p\in D^2_{-},\; x\in\partial U\cap V \}.$$
  $D^2_{-}$ is oriented   by
 the form $ d\theta\wedge dt $, and $\partial U\cap V$ is oriented
 with the orientation of $\partial U$. Hence
 \begin{align}
 \mathcal {T}_{\bf{13}} & =\frac{-i}{2\pi}\int_{A_{\bf{13}}}d\,\log r_{UV}\wedge\big(\omega+d(\alpha(f_t\circ\psi_t))\wedge dt
 \big)^n \\ \notag  &=\frac{-ni}{2\pi}
\int_{A_{\bf{13}}}d\,\log r_{UV}
(f_t\circ\psi_t))\alpha'(\theta)d\theta\wedge dt\wedge\omega^{n-1}
 \\ \notag &=\frac{+ni}{2\pi}\int_{0}^1dt\int_{\partial U\cap
 V}(f_t\circ\psi_t)d\,\log r_{UV}\wedge\omega^{n-1}.
 \end{align}
 In general, if $j<k$
 \begin{equation}\label{1star}
 \mathcal {T}_{\bf{2j-1,2k-1}}=
 \frac{ni}{2\pi}\int_{0}^1dt\int_{A_{jk}}(f_t\circ\psi_t)d\,\log
 r_{jk}\wedge\omega^{n-1},
 \end{equation}
 where $A_{jk}$ is the  set defined in Lemma
 \ref{Lemaux}.

On the other hand
 $$A_{\bf{14}}=\big(\partial U_-\setminus(U_+\cup V_-)\big)\cap
 V_+ =\{[-,p,x]\,|\, p\in D^2_-,\,x\in\partial U\setminus V  \}\cap V_+=\emptyset.$$
 Thus $\mathcal {T}_{\bf{14}}=0$. In general, for $j<k$ the
 integration domain
$A_{\bf{2j-1,2k}}$ is of the form
  $$(\partial
B_{j-}\setminus\cup\,\cdot\,)\cap B_{k+}.$$
 In the union $\cup\,\cdot\,$ appear the sets $B_{k-}$ and $B_{j+}$, hence
$$A_{\bf{2j-1,2k}}\subset (\partial B_{j-}\setminus(B_{j+}\cup
B_{k-}))\cap B_{k+},$$
  and   this set is empty by the same reason
that $A_{\bf{14}}=\emptyset$. Therefore $\mathcal
{T}_{\bf{2j-1,2k}}=0$, for any $j<k$.

The set $A_{\bf{23}}$ is
$$A_{\bf{23}}=( \partial U_+\setminus U_-)\cap V_- =\{ [+,p,x]\,|\, p\in F,\,x\in\partial U\setminus V  \}.$$
As $d\,\log l_{\bf{23}}\wedge\omega^n$ does not contain $d\theta$,
the term $\mathcal {T}_{\bf{23}}$ vanishes. In general, if $j<k$
\begin{align}
\notag A_{\bf{2j,2k-1}}&=
  (\partial
B_{j+}\setminus\cup\,\cdot\,)\cap B_{k-}\subset (\partial
B_{j+}\setminus B_{j-} )\cap B_{k-}  \\ \notag &= \{ [+, p,x]\,|\,
p\in F,\,x\in\partial B_j\setminus B_k  \}.
 \end{align}
  Then
$\mathcal {T}_{\bf{2j,2k-1}}$ vanishes by the same reason that
$\mathcal {T}_{\bf{23}}=0.$

Analogous arguments as the ones explained in the preceding
paragraph show that $\mathcal {T}_{\bf{2j,2k}}=0$, for any $j<k$.

So, apart from the terms $\mathcal {T}_{\bf{ab}}$ considered in
(\ref{2stars}) and in (\ref{1star}), the remainder summands in
(\ref{ITH}) are zero. The theorem follows from (\ref{2stars}) and
(\ref{1star}).

\end{proof}

From the definition of product in $\pi_1(\text{Ham}(M,\omega))$ by
juxtaposition of paths and under the hypotheses of Theorem
\ref{Thefl} is obvious that
$$I:\pi_1(\text{Ham}(M,\omega))\to{\Bbb R}$$
is a group homomorphism.  This fact has been proved in
\cite{L-M-P} for the general case.

\begin{Cor}\label{Corcons} If  $U$ and $V$ are symplectic trivializations of $TM$, with
 $\psi_t(U)=U$, $\psi_t(V)= V$, for all $t$ and $U\cup V=M$ and
$\int_{S^1}(f_t\circ\psi_t) dt$ is a constant $k$  on
 $\partial U\cap V$, then
$$I_{\psi}=J_U\int_U\omega^n+J_V\int_{V\setminus U}\omega^n-nk\langle
c_1(TM)[\omega]^{n-1},\,M \rangle.$$
\end{Cor}

\begin{Cor}\label{Corsim}
If $TM$ is trivial on $U:=M\setminus\{q\}$, where $q$ is a point
of $M$ fixed by $\psi_t$ for all $t$, then
$$I_{\psi}=J_U\int_M\omega^n-n\Big(\int_{S^1}f_t(q) dt\Big)\langle
c_1(TM)[\omega]^{n-1},\,M \rangle.$$
\end{Cor}

 \smallskip

 Now we analyze  the expression for $I_{\psi}$ given in Theorem
  \ref{Thefl} in case of {\it integrable systems}. Let $f$ be the
 normalized Hamiltonian which generates  the loop $\psi$. We
 assume that $(M,\omega,f)$ is completely integrable, with
 $f_1=f,f_2,\dots, f_n$  integrals of motion. We
 suppose that $df_1,\dots, df_n$ are independent at the points
 of $M\setminus P=:V$, where $P$ is a finite union of $2n-2$
 dimensional submanifolds of $M$. We suppose that on $V$ are
 defined action-angle coordinates. We put
 $$Q:=\{x\in P\,|\, \text{dim\,\,Span}\,( df_1(x),\dots,df_n(x))=n-1 \}.$$
 By $Q_1,\dots, Q_k$ are denoted the connected components of $Q$,
 and let $V_j$  be a tubular neighborhood  of $Q_j$ in $M$,
 invariant under $\psi_t$ for all $t$. We assume that on $V_j$ is
 defined a symplectic trivialization of $TM$.
 Then, for each $j$ one can choose a family of tubular
 neighborhoods $\{V_{jb}\subset V_j  \}_{b=1,2\dots}$, such that
 $$\lim_{b\to\infty}\int_{V_{jb}}\omega^n=0.$$
 Lemma  \ref{Lemaux}  applied to the covering $\{V,V_{jb}\}_{j=1,\dots,
 k}$ of $V\cup Q$ gives
$$\langle c_1(M)[\omega]^{n-1},\,[M]\rangle
= \frac{-i}{2\pi}\sum_{j=1}^k\int_{\partial V \cap V_{jb}}d\,\log
r_{VV_{jb}}+\epsilon(b),$$
  where $\epsilon(b)$ goes to $0$ as
$b\to\infty$.

Hence
 \begin{equation}\label{c1int} \langle c_1(M)[\omega]^{n-1},\,[M]\rangle=\sum_{j=1}^k z_j,
 \end{equation}
 with
 \begin{equation}\label{zdef}
 z_j:= \frac{-i}{2\pi}\sum\int_{\partial V \cap
V_{jb}}d\,\log r_{VV_{jb}}\wedge \omega^{n-1}.
\end{equation}

\begin{Prop} Let $(M,\omega,f,f_2,\dots f_n)$ be  an integrable
in which the preceding hypotheses hold, then
$$I_{\psi}=\sum_{j=1}^k z'_j,$$
where $z'_j$ is obtained from the corresponding $z_j$ by inserting
the factor $-nf$ in the integrand of (\ref{zdef}).
\end{Prop}
\begin{proof} The Maslov index $J_V=0$ because of the particular
form of the flow equations in action-angle coordinates. On the
other hand
$$\int_{V_{jb}\setminus (V\cup\dots)}\omega^n=0.$$
 Thus the Proposition follows from Theorem \ref{Thefl}, together
 with (\ref{c1int}) and (\ref{zdef}).
\end{proof}

Similar arguments to the ones involved in this Proposition are
used in Section 4 for  studying the invariant $I$ in Hirzebruch
surfaces.

\smallskip

\section{Examples.}

{\em The invariant $I$ when the manifold is the $2$-sphere.}

 Let $\psi_t$ be the rotation in ${\Bbb R}^3$ around
$\vec{e}_3 $ of angle $2\pi t$ with $t\in[0,\,1]$. Then $\psi_t$
determines a Hamiltonian symplectomorphism of
$(S^2,\omega_{area})$. In fact, the isotopy $\{\psi_t\}$ is
generated by the vector field $\frac{\partial}{\partial\phi}$, and
the function $f$ on $S^2$ defined by
$f(\theta,\phi)=-2\pi\cos\theta=-2\pi z$ is the corresponding
normalized Hamiltonian.

$TS^2$ can be trivialized on $U=D^2_+$, and on $V=D^2_{-}$.
Moreover $\partial  U\cap V$ is the parallel
$\theta=\pi/2+\epsilon$. On $\partial U\cap V$ the function
$f\circ\psi_t$ takes the value $2\pi\sin\epsilon$.
$$\int_U\omega=2\pi(1-k'),\;\;\; \int_{V\setminus U}\omega=2\pi(1+k'),$$
with $k':=\cos (\pi/2+\epsilon)$.

  Furthermore the north pole $n$ and the south
pole $s$ are fixed points of the isotopy $\psi_t$. The rotation
$\psi_t$ transforms the basis $\vec{e}_1,\,\vec{e}_2$ of $T_n S^2$
in
$$(\cos 2\pi t\,\vec{e}_1+\sin 2\pi t\,\vec{e}_2,\;-\sin 2\pi t\,\vec{e}_1+\cos 2\pi t\,\vec{e}_2).$$
So  $J_U$ is  the winding number of the map
$$t\in[0,\,1]\to e^{2\pi ti}\in U(1);$$
That is, $J_U=+1$.

Similarly, by considering the oriented basis
$\vec{e}_2,\,\vec{e}_1$ of $T_s S^2$ it turns out that the Maslov
index $J_V$ of $\psi_t$ is $-1$.

By Corollary \ref{Corcons}
$$I_{\psi}=2\pi(1-k')-2\pi(1+k')-(-2\pi k')\langle
c_1(TS^2),\,S^2\rangle=0.$$

Corollary \ref{Corsim} can also be applied   to determine
$I_{\psi}$. One takes  $U:=S^2\setminus\{s\}$. As
$f(s)=-2\pi(-1)$, we obtain again
 $$I_{\psi}=+4\pi-2\pi\langle
c_1(TS^2),\,S^2\rangle=0.$$

Using formula (\ref{finalfr}) we can determine $I_{\psi}$ again.
Now $V$ is $S^2\setminus\{n,s\}$, $U_1$ is a small polar cap at
$n$ and $U_2$ the symmetric one at $s$. By the symmetry
$$\int_{\partial U_1\cap V}d\,\log r_{U_1V} =
\int_{\partial U_2\cap V}d\,\log r_{U_2V} ,$$
  so $y_1=y_2$. As
$f(n)=-f(s)$, then $I_{\psi}=0$.

This result was expected, because $\pi_1(\text{Ham}(S^2))$ is
isomorphic to ${\Bbb Z}_2$ (see \cite{lP01})  and $I$ is a group
homomorphism.

\vskip .5cm

{\em The invariant $I$ for Hamiltonian loops in ${\Bbb T}^{2n}$.}

We identify the torus ${\Bbb  T}^{2n}$ with ${\Bbb R}^{2n}/{\Bbb
Z}^{2n}$, and we suppose that ${\Bbb  T}^{2n}$ is equipped with
the standard symplectic form $\omega_0$.
 If $\psi_t$
is a Hamiltonian isotopy  of ${\Bbb  T}^{2n}$, it can be written
in the form
$$\psi_t(x^1,\dots,x^{2n})=\big( x^1+\alpha^{1}(t,x^i),\dots, x^{2n}+\alpha^{2n}(t,x^i) \big),$$
where the function $\alpha^j$, for $j=1,\dots, 2n$, is periodic of
period $1$ in each variable: $t,x^1,\dots,x^{2n}$. The vector
fields $\{\frac{\partial}{\partial x^i}\}$ give a symplectic
trivialization of the tangent bundle. In this case the right hand
side of (\ref{subs}) has only one term. The matrix of $(\psi_t)_*$
with respect to $\{\frac{\partial}{\partial x^i}\}$ is
\begin{equation}\label{matxs}
\Big( \delta^j_i+\frac{\partial \alpha^j}{\partial x^i} \Big)\in
Sp(2n,\,{\Bbb R}).
\end{equation}

 First, let us assume that each $\alpha^j$ is a separate  variables
 function; that is, $\alpha^j(t,x^i)=f^j(t) u^j(x^i)$. Since
 $\alpha^1$ takes the same value at symmetric points on opposite
 faces of the cube $I^{2n}$, there is a point $p_1\in I^{2n}$
 such
 $$\frac{\partial u^1}{\partial x^j}(p_1)=0,
 $$
 for all $j$. Hence the first row of the
 matrix (\ref{matxs}) at the point $p_1$ is
 $(1,0,\dots,0)$; that is, the matrix of $(\psi_t)_*(p_1)$ is
 independent of $f^1$ and thus  the Maslov index of $\{(\psi_t)_*(p_1)\}_t$
 does not depend on $f^1$. From (\ref{subs}) it follows that
  $I_{\psi}$ is independent of
 $f^1$. The independence of $I_{\psi}$ with respect to
  $f^i$ is proved in a
 similar way. Thus in order to determine $I_{\psi}$ we can assume that $f^i=0$ for all $i$,
 but in this case
 $I_{\psi}=0$ obviously.

If $\alpha^j$ is sum of two separate variables functions
$$\alpha^j(t,x^i)=f^j(t)u^j(x^i)+g^j(t)v^j(x^i),$$
we take a point $q_1\in I^{2n}$, such that $\frac{\partial
v^1}{\partial x^j}(q_1)=0,$ for all $j$. Then $I_{\psi}$ is
independent of $g^1$. The above reasoning gives $I_{\psi}=0$ in
this case as well.

By the Fourier theory, the original $C^{\infty}$ periodic function
$\alpha^j$ can be approximated (in the uniform $C^{k}$-norm) by a
sum of separated functions of the form $\sum f_a(t)u_a(x^i)$,
where $f_a$ and $u_a$ are $1$-periodic. As $I_{\psi}$ depends only
on the homotopy class of $\psi$, we conclude that $I_{\psi}=0$ for
a general Hamiltonian loop.

\begin{Prop} \label{Itorus}
The invariant $I$ is identically zero on $\pi_1(\text{Ham}({\Bbb
T}^{2n},\,\omega_0))$.
\end{Prop}

This result when $n=1$ is consistent with
 the fact that $\pi_1(\text{Ham}({\Bbb T}^2))=0$ (see \cite{lP01})

\smallskip

{\em Application to  Hirzebruch surfaces.}

 Given $3$
numbers $k,\tau,\mu$, with $k\in{\Bbb Z}_{>0}$, $\tau,\mu\in{\Bbb
R}_{>0}$ and $k\mu<\tau$, the triple $(k,\tau,\mu)$ determine a
Hirzebruch surface $M_{k,\tau,\mu}$ \cite{Au}. This manifold is
the quotient
$$\{z\in {\Bbb C}^4\,\,:\,\, k|z_1|^2+|z_2|^2+|z_4|^2=\tau/\pi,\;
|z_1|^2+|z_3|^2=\mu/\pi\}/{\Bbb T}^2,$$ where the ${\Bbb
T}^2$-action is given by
$$(a,b)\cdot(z_1,z_2,z_3,z_4)=(a^kbz_1,az_2,bz_3,az_4),$$
for $(a,b)\in {\Bbb T}^2$. The map
$$[z_1,z_2,z_3,z_4]\mapsto ([z_2:z_4],[z_2^kz_3:z_4^kz_3:z_1])$$
allows us to represent $M_{k,\tau,\mu}$ as a submanifold of ${\Bbb
CP}^1\times {\Bbb CP}^2$. On the other hand the usual symplectic
structures on
 ${\Bbb CP}^1$ and ${\Bbb CP}^2$ induce a symplectic form $\omega$ on
 $M_{k,\tau,\mu}$, and the following ${\Bbb T}^2$-action on ${\Bbb CP}^1\times
{\Bbb CP}^2$
$$(a,b)([u_0:u_1],[x_0:x_1:x_2])=([au_0:u_1],[a^kx_0:x_1:bx_2])$$
gives rise to a toric structure on $M_{k,\tau,\mu}$. In terms of
the Delzant construction $(M_{k,\tau,\mu},\,\omega)$ is associated
to the trapezoid in $({\Bbb R}^2)^*$ whose not oblique edges are
$\tau,\mu, \lambda:=\tau-k\mu$, \cite{Gui} (see Figure 2).
Moreover $\lambda$ is the value that the symplectic form $\omega$
takes on the exceptional divisor, $\{[z]\in M \,|\, z_3=0\}$, of
$M:=M_{k,\tau,\mu}$. And $\omega$ takes the value $\mu$ on the
class of the fibre in the fibration $M\to {\Bbb C}P^1$.

\begin{figure}[htbp]
\begin{center}
\epsfig{file=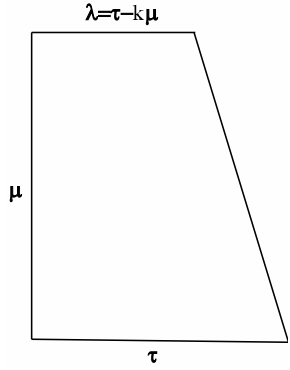,height=5cm}
\end{center}
\caption[Figure 2]{\small Delzant polytope associated to $M $. }
\end{figure}

Since $M$ is a toric manifold, the ${\Bbb T}^2$-action define
symplectomorphisms of $M$. More precisely, let $\psi_t$ the
diffeomorphism of $M$ defined by
\begin{equation}\label{psi}
\psi_t[z_1,z_2,z_3,z_4]=[z_1e^{2\pi it},z_2,z_3,z_4].
\end{equation}
$\psi=\{\psi_t\,:\,t\in [0,1]\}$ is a loop of Hamiltonian
symplectomorphisms of $(M,\omega)$. Similarly we have
\begin{equation}\label{psitil}
\tilde\psi_t[z_1,z_2,z_3,z_4]=[z_1,z_2e^{2\pi it},z_3,z_4],
\end{equation}
 and the corresponding loop $\tilde\psi$ in
$\text{Ham}(M,\omega).$

Using
 Theorem \ref{Thefl} we shall calculate
the values of $I_{\psi}$ and $I_{\tilde\psi}$ in terms of
$\lambda,$ $\tau$ and $k$. The result is stated in Theorem
\ref{ThemHiz} below. The most laborious point in the proof of the
following Theorem is to obtain Darboux charts for $M$ which give
rise to simple transition functions for $\text{det}(TM)$.

\begin{Thm}\label{ThemHiz}
 Let $\psi$ and $\tilde\psi$ be the loops of
symplectomorphisms of the Hirzebruch surface
$(M_{k,\tau,\mu},\omega)$, defined by (\ref{psi}) and
(\ref{psitil}) respectively, then
$$I_{\psi}=\frac{2k\mu^2}{3}\Big(1-\frac{\mu}{2\lambda+k\mu}\Big),\;\;\text{and}\;\;
I_{\tilde\psi}=\frac{-k^2\mu^2}{3}\Big(1-\frac{\mu}{2\lambda+k\mu}\Big).$$
$\lambda$ being $\tau-k\mu$.
\end{Thm}

\begin{proof}

We will define a Darboux atlas on $M$. First we consider the
following covering for $M$
$$U_1=\{[z]\in M\,:\, z_3\ne 0\ne z_4  \},\;\; U_2=\{[z]\in M\,:\, z_1\ne 0\ne z_4
\}$$
$$U_3=\{[z]\in M\,:\, z_1\ne 0\ne z_2  \},\;\; U_4=\{[z]\in M\,:\, z_2\ne 0\ne
z_3 \}.$$
 We set $z_j=\rho_je^{i\theta_j}$, with $\rho_j=|z_j|$, and on
 $U_1$ introduce the coordinates $(x_1,y_1,a_1,b_1)$ by the
 formulae
 $$x_1+iy_1=\rho_1e^{i\varphi_1},\;\;
 a_1+ib_1=\rho_2e^{i\varphi_2},\;\;
\varphi_1=\theta_1-\theta_3-k\theta_4,\;\;
\varphi_2=\theta_2-\theta_4.$$ Then $\omega$ on $U_1$ can be
written $\omega=dx_1\wedge dy_1 +da_1\wedge db_1.$

On $U_2$ we consider the Darboux coordinates $(x_2,y_2,a_2,b_2)$,
with
$$x_2+iy_2=\rho_3e^{i\xi_3},\;\;
 a_2+ib_2=\rho_2e^{i\xi_2},\;\;
\xi_2=\theta_2-\theta_4,\;\; \xi_3=\theta_3-\theta_1+k\theta_4.$$

On $U_3$ we put
$$x_3+iy_3=\rho_3e^{i\chi_3},\;\;
 a_3+ib_3=\rho_4e^{i\chi_4},\;\;
\chi_3=\theta_3-\theta_1+k\theta_2,\;\;
\chi_4=\theta_4-\theta_2,$$ and $\omega=dx_3\wedge dy_3
+da_3\wedge db_3.$

Finally, on $U_4$ we set
$$x_4+iy_4=\rho_1e^{i\zeta_1},\;\;
 a_4+ib_4=\rho_4e^{i\zeta_4},\;\;
\zeta_1=\theta_1-\theta_3-k\theta_2,\;\;
\zeta_4=\theta_4-\theta_2,$$ and $\omega=dx_4\wedge dy_4
+da_4\wedge db_4.$

The  normalized Hamiltonian function for $\psi_t$ is
$f=\pi\rho_1^2-\kappa$, where $\kappa$ is a constant determined by
the condition $\int_M f\omega^2=0$. Straightforward  calculations
give
$$\int_M\omega^2=\mu(2\tau-k\mu),\;\;\text{and}\;\;
\int_M\pi\rho_1^2\omega^2=\frac{\mu^2}{3}(3\tau-2k\mu).$$ So
 \begin{equation}\label{kapp}
 \kappa=\frac{\mu}{3}\Big(\frac{3\lambda+k\mu}{2\lambda+k\mu}\Big).
 \end{equation}

It is not easy to determine the transition function of
$\text{det}(TM)$ that corresponds to the coordinate transformation
$(x_i,y_i,a_i,b_i)\to (x_j,y_j,a_j,b_j)$; that is why we will
introduce polar coordinate on subsets of the domains $U_j$.

Given $0<\epsilon<<1$, for $j=1,2,3,4$ we put
$$B_j=\{[z]\in
U_j\,:\, |z_j|<2\epsilon \}\;\;\text{and}\;\; B_0=\{[z]\in M\,:\,
|z_j|>\epsilon \;\text{for all}\; j\}.$$ On $B_0$ are well-defined
the coordinates
$(\frac{\rho_1^2}{2},\varphi_1,\frac{\rho_2^2}{2},\varphi_2)$, and
in this coordinates
$$\omega= d\Big(\frac{\rho_1^2}{2}\Big)\wedge d\varphi_1+
d\Big(\frac{\rho_2^2}{2}\Big)\wedge d\varphi_2.$$

On $B_j$ ($j=1,2,3,4$) we consider the Darboux coordinates
$(x_j,y_j,a_j,b_j)$ defined above. Then $B_0,B_1,B_2,B_3,B_4$ is a
Darboux atlas for $M$.
 We
assume that $M$ is endowed with the orientation given by
$\omega^2$. This orientation agrees  on $B_0$ with the one defined
by $d\rho_1^2\wedge d\varphi_1\wedge d\rho_2^2\wedge d\varphi_2$.

It is evident that $\psi_t(B_i)=B_i$, for $i=0,1,2,3,4$. Since
$\psi_t$ on $B_0$ is simply the translation
$\varphi_1\to\varphi_1+2\pi t$ of the variable $\varphi_1$, the
Maslov index $J_0$ of $\psi$ in the trivialization defined on
$B_0$ vanishes.

As $B_j$ (for $j=1,2,3,4$)  has "infinitesimal size" and $J_0=0$,
the expression for $I_{\psi}$ of Theorem \ref{Thefl} can be
written
\begin{equation}\label{HizI}
I_{\psi}=\sum_{i<k}N_{ik}+O(\epsilon)
\end{equation}
Since $I_{\psi}$ is obviously independent of the coordinates, it
follows from (\ref{HizI}) that $N_{ik}$ is independent, up to
order $\epsilon$, of the chosen Darboux coordinates in $B_j$, for
$j=1,2,3,4$. Moreover $N_{ik}$ with $0\ne i<k$ is also of order
$\epsilon$.

On the other hand, if we substitute $B_j$ by
$$B'_j=\{[z]\in
B_j\,:\, |z_r|>\epsilon,\; r\ne j\}$$
 in the definition of
$N_{ik}$  (see Theorem \ref{Thefl}) the new $N_{ik}$ differs from
the old one in a quantity of order $\epsilon$. As on $B'_1$ the
variable $\rho_2\ne 0$, we can consider the Darboux coordinates
$$(x_1,y_1,\frac{\rho_2^2}{2},\varphi_2)$$
 on $B'_1$.
Since $\rho_3\ne 0$ on $B'_2$ we take the coordinates
$(a_2,b_2,\frac{\rho_3^2}{2},\xi_3)$ on $B'_2$. Similarly we will
adopt the following coordinates:
$(x_3,y_3,\frac{\rho_4^2}{2},\chi_4)$ on $B'_3$ and
$(\frac{\rho_1^2}{2},\zeta_1,a_4,b_4)$ on $B'_4$.

Taking into account the preceding arguments
\begin{equation}\label{Ipsia}
I_{\psi}=\sum_{j=1}^4 N'_{0j}+O(\epsilon),
\end{equation}
where
\begin{equation}\label{N0j}
N'_{0j}=\frac{i}{\pi}\int_{A'_{0j}}f d\,\log r_{0j}\wedge\omega
\end{equation}
and
$$A'_{0j}=\{[z]\in M\,:\, |z_r|>\epsilon, \;\text{for all}\;
r\ne j \;\text{and}\; |z_j|=\epsilon \}.$$
 The submanifold
$A'_{0j}$ is oriented as a subset of $\partial B_0$; that is, with
the orientation induced by the one of $B_0$.

Next we determine the value of $N'_{01}$. To know the transition
function $r_{01}$ one needs the Jacobian matrix $R$ of the
transformation
$$(x_1,y_1,\frac{\rho_2^2}{2},\varphi_2)\to
(\frac{\rho_1^2}{2},\varphi_1,\frac{\rho_2^2}{2},\varphi_2)$$ in
the points of $A'_{01}$; with $\rho_1^2=x_1^2+y_1^2$,
$\varphi_1=\tan^{-1}(y_1/x_1)$. The non trivial block of $R$ is
the diagonal one
$$\begin{pmatrix} x_1& y_1 \\
r& s \\
\end{pmatrix},$$
 with $r={-y_1}{(x_1^2+y_1^2)}^{-1}$ and
$s={x_1}{(x_1^2+y_1^2)}^{-1}$. The non-real eigenvalues of $R$ are
$$\lambda_{\pm}=\frac{x_1+s}{2}{\pm}\frac{i\sqrt{4-(s+x_1)^2}}{2}.$$
On $A'_{01}$  these non-real eigenvalues occur when $(s+x_1)^2<2$,
that is, if $|\cos
\varphi_1|<2\epsilon(\epsilon^2+1)^{-1}=:\delta$. If $y_1>0$ then
$\lambda_{-}$ of the first kind (see \cite{S-Z}) and $\lambda_+$
is of the first kind, if $y_1<0$.

Hence, on $A'_{01}$,
\begin{equation}\notag \rho(R)=\begin{cases}
\lambda_+|\lambda_+|^{-1}=x+iy,
&\text{if $\,|\cos\varphi_1|<\delta\,$ and $\,y_1<0$;}\\
\lambda_-|\lambda_-|^{-1}=x-iy,
&\text{if $\,|\cos\varphi_1|<\delta\,$ and $\,y_1>0$;}\\
\pm 1,&\text{otherwise.}
\end{cases}
\end{equation}
where $x=\delta^{-1}\cos\varphi_1$, and $y=\sqrt{1-x^2}$.

If we put $\rho(R)=e^{i\gamma}$, then $\cos\gamma=\delta^{-1}\cos
\varphi_1$ (when $|\cos\varphi_1|<\delta$), and
\begin{equation}\notag
\sin\gamma=
\begin{cases}-\sqrt{1-\cos^2\gamma},&\text{if
$\,\sin\varphi_1>0$;}\\
\sqrt{1-\cos^2\gamma},&\text{if $\,\sin\varphi_1<0$.}
\end{cases}
\end{equation}
So when $\varphi_1$ runs anticlockwise from $0$ to $2\pi$,
$\gamma$ goes round clockwise the circumference; that is,
$\gamma=h(\varphi_1)$, where $h$ is a function such that
\begin{equation}\label{h1}
h(0)=2\pi,\;\; \text{and}\;\; h(2\pi)=0.
 \end{equation}
 As $r_{01}=\rho(R)$,
then $d\,log \,r_{01}=id h$.

On $A'_{10}$ the form   $\omega$ reduces to $(1/2)d\rho_2^2\wedge
d\varphi_2$. From (\ref{N0j}) one deduces
\begin{equation}\label{N01a}
N'_{01}=\frac{i}{2\pi}\int_{A'_{01}}i f dh\wedge d\rho_2^2\wedge
d\varphi_2.
\end{equation}

 On the other hand according to the convention about orientations,
$\{[z]\,:\, |z_1|=\epsilon\}$ as subset of $\partial B_0$ is
oriented by $-d\varphi_1\wedge d\rho_2^2\wedge d\varphi_2.$ And on
$A'_{01}$ the Hamiltonian function $f=-\kappa+O(\epsilon)$. Then
it follows from (\ref{N01a}) together with (\ref{h1})
\begin{equation}\label{N01}
N'_{01}=2\tau\kappa+O(\epsilon).
\end{equation}

The contributions $N'_{02}, N'_{03}, N'_{04}$ to $I_{\psi}$ can be
calculated in a similar way. One obtains the following results up
to addends of order $\epsilon$:
\begin{equation}\label{contr}
N'_{02}=2\mu\kappa
-\mu^2,\;\;N'_{03}=2\lambda(\kappa-\mu),\;\;N'_{04}=\mu(2\kappa-\mu).
\end{equation}
As $I_{\psi}$ is independent of $\epsilon$, it follows from
(\ref{Ipsia}), (\ref{N01}), (\ref{contr}) and (\ref{kapp})
$$I_{\psi}=\frac{2k\mu^2}{3}\Big(1-\frac{\mu}{2\lambda+k\mu}
\Big).$$

Next we consider the loop $\tilde\psi$; the corresponding
normalized Hamiltonian function is $\tilde
f=\pi\rho_2^2-\tilde\kappa$, where
\begin{equation}\label{kapptil}
\tilde\kappa=\frac{3\lambda^2+3k\lambda\mu+k^2\mu^2}{3(2\lambda+k\mu)}.
\end{equation}
As in the preceding case
\begin{equation}\label{Ipsitil}
I_{\tilde\psi}=\sum_{j=1}^4 \tilde N'_{0j} +O(\epsilon),
\end{equation}
where
$$\tilde N'_{0j}=\frac{i}{\pi}\int_{A'_{0j}}\tilde f d\,\log
r_{0j}\wedge\omega.$$

The expression for $\tilde N'_{01}$ can be obtained from
  (\ref{N01a}) substituting  $f$ for $\tilde f$; so
  \begin{equation}\label{N01t}
  \tilde N'_{01}=\tau(2\tilde\kappa-\tau)+O(\epsilon).
  \end{equation}
  Analogous calculations give the following values for the $\tilde N'_{0j}$'s, up to summands of order $\epsilon$
\begin{equation}\label{Ntildes}
\tilde N'_{02}=2\mu\tilde\kappa,\;\; \; \tilde
N'_{03}=\lambda(2\tilde\kappa-\lambda),\;\;\; \tilde
N'_{04}=\mu(2\tilde\kappa-k\mu-2\lambda).
 \end{equation}

  From (\ref{kapptil}), (\ref{Ipsitil}), (\ref{N01t}) and (\ref{Ntildes}) it
  follows the value for $I_{\tilde\psi}$ given in the statement of Theorem.

\end{proof}

{\it Remark.}  In \cite{A-Mc} is proved that
$\pi_1(\text{Ham}(M))={\Bbb Z}$ when $k=1$, therefore the quotient
of $I_{\psi}$ by $I_{\psi'}$, for arbitrary Hamiltonian loops of
symplectomorphisms, is a rational number. For the particular loops
considered in Theorem \ref{ThemHiz} the quotient
$I_{\tilde\psi}/I_{\psi}$ equals $-k/2$, so Theorem \ref{ThemHiz}
is consistent with the result of Abreu and McDuff.


\end{document}